\documentclass[twocolumn]{svjour3}          
\smartqed  
\usepackage{graphicx}

\usepackage{latexsym,bm,amssymb,amsfonts,amsbsy,amsmath,graphicx,color,marvosym}
\usepackage{pifont}
\usepackage{url}
\usepackage{amsmath}
\usepackage{amsfonts}
\usepackage{amssymb}
\usepackage{fancyhdr}
\usepackage{pstricks,pst-node,pst-text,pst-3d}
\def\go{\hbox{$\displaystyle{\mbox{\ding{172}}}$}}
\def\gob{\hbox{$\displaystyle{\mbox{\ding{182}}}$}}
\def\G1{\hbox{$\displaystyle{\mbox{\ding{172}}}$}}

\def \RR {{\mathbb{R}}}

\def \FF {{\mathbb{F}}}

\def\eps{\varepsilon}

\newrgbcolor{marrone}{0.6156 0.3164 0.0313}
\newrgbcolor{ts}{0.50 0.598 0.781}
\newrgbcolor{verde}{0 0.4 0}
\newrgbcolor{LBlue}{0.58 0.65 0.9}
\newrgbcolor{Blue}{0.1 0.00 0.60}
\newrgbcolor{rosso}{0.9633    0.1641    0.0039}

\journalname{~}
\begin{document}

\title{On the use of the Infinity Computer architecture to set up a dynamic precision floating-point arithmetic
}


\author{Pierluigi Amodio \and Luigi Brugnano  \and Felice Iavernaro$^{\textrm{\Letter}}$ \and Francesca Mazzia
}


\institute{P. Amodio \at
              Dipartimento di Matematica, Universit\`a di Bari, Italy \\
              \email{pierluigi.amodio@uniba.it}           
           \and
           L. Brugnano \at
              Dipartimento di Matematica e Informatica ``U. Dini'', Universit\`a di Firenze, Italy \\
              \email{luigi.brugnano@unifi.it}
           \and
           F. Iavernaro \at
              Dipartimento di Matematica, Universit\`a di Bari, Italy \\
              \email{felice.iavernaro@uniba.it} 
           \and
           F. Mazzia \at
           Dipartimento di Informatica, Universit\`a di Bari, Italy \\
           \email{francesca.mazzia@uniba.it} 
}

\date{}

\maketitle

\begin{abstract}
We devise a variable precision floating-point arithmetic by exploiting the framework provided by the Infinity Computer. This is a computational platform implementing the Infinity Arithmetic system, a positional numeral system 
which can handle both infinite and infinitesimal quantities symbolized by the positive and negative finite powers of the radix $\go$. The computational features offered by the Infinity Computer allows us to dynamically change the accuracy of representation and floating-point operations during the flow of a computation. When suitably implemented, this possibility turns out to be particularly advantageous when solving ill-conditioned problems. In fact, compared with a standard multi-precision arithmetic, here the accuracy is improved only when needed, thus not affecting that much the overall computational effort. An illustrative example about the solution of a nonlinear equation is also presented.

\end{abstract}

\section{Introduction}
\label{intro}
The Arithmetic of Infinity was introduced by Y.D. Ser\-ge\-yev with the aim of devising a new coherent computational environment able to handle finite, infinite and infinitesimal quantities, and to execute arithmetical operations with them. It is based on a positional numeral system with the infinite radix   $\go$, called \textit{grossone} and representing, by definition,  the  number of elements of the set of natural numbers $\mathbb{N}$ (see, for example, \cite{Ser08,Ser09} and the survey paper \cite{EMS}). Similar to the standard positional notation for finite real numbers, a number in this  system is recorded as
$$
c_{{m}}
\G1^{p_{m}}    \ldots   c_{{1}} \G1^{p_{1}} c_{{0}} \G1^{p_{0}}
c_{{-1}} \G1^{p_{-1}} \ldots c_{{-k}}
\G1^{p_{-k}},
$$
with the obvious meaning
\begin{equation}
c_{{m}}
\G1^{p_{m}}  +  \ldots + c_{{1}} \G1^{p_{1}} +c_{{0}} \G1^{p_{0}}  + c_{{-1}} \G1^{p_{-1}} + \ldots   + c_{{-k}}
\G1^{p_{-k}}.
\label{numberc}
\end{equation}
The coefficients $c_i$, called \textit{grossdigits}, are real numbers while the \textit{grosspowers} $p_i$, sorted in decreasing order
$$
p_{m} >  \ldots    > p_{1} > p_0=0 > p_{-1}  > \ldots >   p_{-k},
$$
may be finite, infinite or infinitesimal even though, for our purposes, only finite integer grosspowers will be considered. 

Notice that, since $\go^0=1$ by definition, the set of real numbers and the related operations are naturally included in this new system. In this respect, the Arithmetic of Infinity should be perceived as a more powerful tool that improves the ability of observing and describing mathematical outcomes that the standard numeral system could not properly handle. In particular, the new system allows us to better inspect the nature of the infinite objects we are dealing with. For example, while $\infty+1=\infty$ in the standard thinking, if we are in the position to specify as, say $\go$, the kind of infinity we are observing  using the new methodology, such an equality could be better replaced with $\go+1>\go$. According to the principle that \textit{the part is less than the whole}, this novel perception of the infinite dimensionality has proved successful in resolving a number of paradoxes involving  infinities and infinitesimals,  the most famous being Hilbert's paradox of the Grand Hotel (see \cite{Ser08,Ser09}).

The Arithmetic of Infinity paradigm is rooted in three methodological postulates and its consistency has been rigorously recognized in  \cite{Lo15}. Its theoretical and practical implications are formidable also considering that the final goal is to make the new computing system available through a dedicated processing unit. The computational device that implements the Infinity Arithmetic has been called \textit{Infinity Computer} and is patented in EU, USA, and Russia (see, for example, \cite{Sergeyev_patent}).

Among the many fields of research this new methodology has been successfully applied, we  mention {\em numerical differentiation and optimization}  \cite{DeLeone,Num_dif,Zilinskas},  {\em numerical solution of differential equations} \cite{Se13,AmIaMaMuSe16,SeMuMaIaAm16,MaSeIaAmMu16,IaMaMuSe19}, {\em models for percolation and biological processes}  \cite{DeBartolo,Iudin},  {\em cellular automata} \cite{Iudin_2,DAlotto_2}.\footnote{For further references and applications see  the  survey \cite{EMS}.}

The aim of the present study is to devise a \textit{dynamic precision} floating-point arithmetic by exploiting the computational platform provided by the Infinity Computer. In contrast with standard variable precision arithmetics, here not only may the accuracy be dynamically changed during the execution of a given algorithm, but variables stored with different accuracies may be combined through the usual algebraic operations. This strategy is explored and addressed to the accurate solution of ill-conditioned/unstable problems \cite{BrMaTr11,IaMaTr06}.

One interesting application is the possibility of handling ill-conditioned problems or even of implementing algorithms which are labeled as unstable in standard floating-point arithmetic.\footnote{First results on handling ill-conditioning using the Infinity Computer may be  found in \cite{GaGiMu18,SeKvMu18}.} One example in this direction has been illustrated in \cite{AmBrIaMa20}. It consists in the use of the iterative refinement to improve the accuracy of a computed solution to an ill-conditioned linear system until a prescribed input accuracy is achieved. 

The paper is organized as follows. In the next section we highlight those features of the Infinity Computer that play a key role to set up the variable-precision arithmetic. This latter is discussed in Section \ref{sec:3} together with a few illustrative examples. 
As an application in Numerical Analysis, in Section \ref{sec:5} we consider the problem of finding the zero of a nonlinear function affected by ill-conditioning issues. Finally, some conclusions are drawn in Section \ref{sec:6}.

\section{Background}
\label{sec:2}
As is the case with the standard floating-point arithmetic,  the Infinity Computer handles both numbers and operations numerically (not symbolically). Consequently, it is prone to efficiently afford the massive amount of computation needed while solving a wide variety of real-life problems. On the other hand,  a roundoff error proportional to the machine accuracy is generated during representation of data (i.e. the coefficients $c_i$ and $p_i$ in (\ref{numberc})) and execution of the basic operations. We will give a more detailed description about how the representation of grossdigits and the floating-pont operations should be carried out in the next section. Here, for sake of simplicity, we will neglect these sources of errors. 

The grossnumbers that will be considered in the sequel are those that admit an expansion in terms of integer powers of $\go^{-1}$ and, thus, take the form
\begin{equation}
\label{X}
X \,=\,\sum_{j=0}^{T} c_j\go^{-j},
\end{equation}
where $T$ denotes the maximum order of infinitesimal appearing in $X$.
 For this special set, the arithmetic operations on the Infinity Computer follow the same rules defined for the polynomial ring. For example, given the two grossnumbers 
\begin{equation}
\label{XY}
X=x_0\go^0+x_1\go^{-1}, ~~ Y=y_0\go^0+y_1\go^{-1}+y_2\go^{-2},
\end{equation}
we get
\begin{equation*}
\label{add}
X+Y=(x_0+y_0)\go^0+(x_1+y_1)\go^{-1}+y_2\go^{-2},
\end{equation*}
\begin{equation*}
\label{prod}
\begin{array}{rcl}
X\cdot Y&=&x_0y_0\go^0+(x_0y_1+x_1y_0)\go^{-1}\\&&+(x_0y_2+x_1y_1)\go^{-2}+x_1y_2\go^{-3},
\end{array}
\end{equation*}
and analogously for the division $X/Y$.
Notice that, on the Infinity Computer, variables may coexist  with  different storage requirements. Taking aside the (negative) powers of $\go$ that, as we will see, need not to be stored in our usage, the variable $Y$ displays infinitesimals quantities up to the order $2$, thus requiring one extra record to store the grossdigit $y_2$, if compared with the variable $X$ that only contains a first order infinitesimal. This circumstance also influences the computational complexity associated with each single floating-point operation.  As a consequence of the different amount of memory allocated for storing grossnumbers,  the global computational complexity associated  with  a given algorithm performed on the Infinity Computer, cannot be merely estimated in terms of how many flops are executed, but should also take into account how many grossdigits are involved in each operation.    

If $X$ is chosen as in (\ref{X}), we denote by $X^{(q)}$ its section obtained by neglecting, in the sum,  all the infinitesimals of order greater than $q$, that is  
\begin{equation}
\label{secX}
X^{(q)}= \sum_{j=0}^q c_j \go^{-j}.
\end{equation}
For example, choosing $q=0$ and $X$ and $Y$ as in  (\ref{XY}), we see that $X^{(0)}+Y^{(0)} =  x_0+y_0$ and $X^{(0)}\cdot Y^{(0)} =  x_0y_0$ would resemble the floating-point addition and multiplication in standard arithmetic, respectively, while additional effort is needed if other powers of $\go^{-1}$ are successively involved. More precisely,
the computational cost associated with a single operation of two grossnumbers will depend on how many infinitesimal are considered. Assuming $q<p$ and denoting by $d_j$ the grossdigits associated with $Y$, for the two sections $X^{(q)}$ and $Y^{(p)}$, the addition 
\begin{equation}
\label{Xq+Yq}
X^{(q)}+Y^{(p)} = \sum_{j=0}^q (c_j+d_j) \go^{-j} + \sum_{j=q+1}^p d_j \go^{-j}
\end{equation}
requires $q+1$ additions of grossdigits, while the multiplication 
\begin{equation}
\label{Xq*Yq}
\begin{array}{rcl}
X^{(q)}\cdot Y^{(p)} &=& \displaystyle \sum_{j=0}^q \sum_{i=0}^j c_id_{j-i} \go^{-j} +
 \sum_{j=q+1}^p \sum_{i=0}^q c_id_{j-i} \go^{-j} \\[.5cm]
&&\displaystyle  + \sum_{j=p+1}^{p+q} ~ \sum_{i=j-p}^q c_id_{j-i} \go^{-j} \\[.5cm]
&=& \displaystyle \sum_{j=0}^{q+p} ~~ \sum_{i=\max\{0,j-p\}}^{\min\{q,j\}} c_id_{j-i} \go^{-j},
\end{array}
\end{equation} 
amounts to $(q+1)(p+1)$ multiplications and $qp-q(q-1)/2$ additions/subtractions of grossdigits.\footnote{The division algorithm is described in Section \ref{sec:3} and therefore is not discussed here.} It is worth noticing that, since in both operations all the coefficients of $\go^{-j}$ may be independently calculated, there is room for a huge parallelization. We will not consider this aspect in detail in the present study.

\section{A variable-precision representation of floating-point numbers on the Infinity Computer}
\label{sec:3}
Grossnumbers of the form (\ref{X}) and their sections (\ref{secX}) form the basis of the new floating-point arithmetic where numbers with a different accuracy may be simultaneously represented and combined. The idea is to let $\go^{-1}$ and its powers act as  \textit{machine infinitesimal quantities} when related to the classical floating-point system. These infinitesimal entities, if suitably activated or deactivated, may be conveniently exploited to increase or decrease the required accuracy during the flow of a given computation. This strategy may be used to automatically detect ill-conditioning issues during the execution of a code that solves a given problem, and to change the accuracy accordingly, in order to optimize the overall computational effort under the constrain that the resulting error in the output solution should fit a given input tolerance. A formal introduction of the new dynamic precision arithmethic is discussed hereafter. 
\subsection{Machine numbers and their storage in the Infinity Computer}
Let $t$ and $T$ be two given non-negative integers and $N=(T+1)(t+1)-1$. The set of machine numbers we are interested in is given by 
\begin{equation}
\label{FF}
\FF = \displaystyle  \left\{ X\in \RR ~\big|~ X=\pm \beta^p\sum_{i=0}^N d_i \beta^{-i} \right\} \cup \{0\},
\end{equation}
where $\beta\ge 2$ denotes the base of the numeral system, the integer $p$ is the exponent ranging in a given finite interval, and $d_i$ are the significant digits, with $d_0\not =0$ (normalization condition). Starting from  $d_0$, we group the digits $d_i$ in $T+1$ adjacent strings each of length $t+1$:
\begin{equation}
\label{Xmt}
\begin{array}{rcl}
X &=& \pm \beta^p \, \underbrace{d_0.d_1\cdots d_t}_{t+1}\underbrace{d_{t+1}\cdots d_{2t+1}}_{t+1}\cdots \\ && 
\underbrace{d_{j(t+1)}\cdots d_{(j+1)(t+1)-1}}_{t+1} \cdots \underbrace{d_{T(t+1)}\cdots d_{(T+1)(t+1)-1}}_{t+1}\\
&=& \displaystyle \pm \beta^p \sum_{j=0}^T \beta^{-j(t+1)} \sum_{i=0}^t d_{j(t+1)+i} \beta^{-i}.
\end{array}
\end{equation}
The representation of the numbers $X$ as in (\ref{FF}), under the shape (\ref{Xmt}), suggests an interesting application of the Infinity Computer. Introducing the new symbol $\gob$, called \textit{dark grossone}, as
\begin{equation}
\label{gob}
\gob=\beta^{t+1},
\end{equation}
and setting  
\begin{equation}
\label{cj}
c_j=\sum_{i=0}^t d_{j(t+1)+i} \beta^{-i},
\end{equation}
the number $X$ in (\ref{Xmt}) may be rewritten as 
\begin{equation}
\label{Xm}
X = \displaystyle \pm \beta^p \sum_{j=0}^T c_j\gob^{-j}.
\end{equation}
Its section of index $q$ is then given by
\begin{equation}
\label{secXm}
X^{(q)}=\displaystyle \pm \beta^p \sum_{j=0}^q c_j\gob^{-j}.
\end{equation}
We assume that a real number $x$ is represented by a floating-point number $X$ in the form (\ref{Xm}) by truncating or rounding it to the nearest even, after the digit $d_N$. This is the most attainable accuracy during the data representation phase but, in general, a lower accuracy (and hence faster execution times) will be required while processing the data, which will be achieved by involving sections of $X$ of suitable indices $q$ during the computations.   

Echoing the symbol $\go$,  the new symbol $\gob$ emphasizes the formal analogy between a machine number and a grossnumber (compare (\ref{Xm}) with (\ref{X}) and (\ref{secXm}) with (\ref{secX})). This correspondence suggests that the computational platform provided by the Infinity Computer may be conveniently exploited to host the set $\FF$ defined at (\ref{FF}) and to execute operations on its elements using a novel methodology. This is accomplished by formally identifying the two symbols, which means that, though they refer to two different definitions, they are treated in the same way in relation to the storage and execution of the basic operations. In accord with the features  outlined in Section \ref{sec:2}, the Infinity Computer will then be able to:
\begin{itemize}
\item[(a)] store floating-point numbers at different accuracy levels, by involving  different infinitesimal terms,  
according to the need;
\item[(b)] easily access to sections of floating-point numbers as defined in (\ref{secXm});
\item[(c)] perform computations involving numbers stored with different accuracies.
\end{itemize} 
The affinity between the meaning of the two symbols goes even beyond what has been stated above.  We have already observed that the case $q=0$ in (\ref{secXm}) resembles the standard set of floating-point numbers with $t+1$ significant figures. This means that when the Infinity Computer works with numbers of the form $X^{(0)}$ it precisely matches the structure designed following the principles of the IEEE 754 standard. In this mode, the operational accuracy is set at its minimum value and the upper bound on the relative error due to rounding (unit roundoff) is $\gob^{-1}$. In other words, $\gob^{-1}$ will be perceived as an infinitesimal entity which cannot be handled unless we let numbers in the form $X^{(1)}$ come into play. This argument  can then be iterated to involve $\gob^{-i}$, $i=2,\dots,T$.  Mimicking the same concept expressed by the use of $\go$, negative powers of $\gob$ act like lenses to observe and combine numbers using different accuracy levels. 

\begin{remark}
What about the role of $\gob$ as an infinite-like quantity?  Consider again the basic operational mode with numbers in the form $X^{(0)}$. If we ask the computer to count integer numbers according to the scheme
\begin{verbatim}
n=0
while n+1>n
   n=n+1
end
\end{verbatim}
it would stop at $n=\gob$, yielding a further similarity with the definition of $\go$ in the Arithmetic of Infinity. Again, involving sections of higher index, the counting process could be safely continued. 
\end{remark}
In conclusion, the role of $\gob$ could be interpreted as an inherent feature of the machine architecture which, consistently with the Infinity Arithmetic methodology, could activate suitable powers of $\gob$ to get, when needed,  a better perception of numbers. The examples included in the sequel further elucidate this aspect. 

\subsection{Floating-point operations}
We have seen that, through the formal identification of $\gob$ with $\go$, it is possible to store the elements of $\FF$ as if they were  grossnumbers and, consequently, to take advantages of the facilities provided by the Infinity Computer in accessing their sections and  performing the four basic operations on them, according to the rules described in Section \ref{sec:2} (see, for example, (\ref{Xq+Yq}) and (\ref{Xq*Yq})). For these reasons, in the sequel, we shall use $\go$ in place of $\gob$ when working on the Infinity Computer, even though,  due to the finite nature of $\gob$, the result of a given operation may not be in the form (\ref{secXm}), so that a normalization procedure has to be considered. Hereafter, we report a few examples in order to elucidate this aspect. For all cases, a binary base has been adopted for data representation.
\paragraph{Addition.} Set $t=3$  and $T=2$ (three grossdigits each with four significant digits), and   
consider the sum of the two floating-point normalized numbers:
$$
\begin{array}{rcl}
X &=&  2^0 \cdot 1.11010101110,  \\
Y &=&  2^{-3} \cdot  1.11111001011. 
\end{array}
$$
Table \ref{tab1} summarizes the procedure by a sequence of commented steps.
\begin{table*}
\caption{Scheme of the addition of two positive floating-point numbers.}
$$
\begin{array}{l@{\quad}|l@{\quad}|l@{\quad}|l@{\quad}|l@{\quad}|l@{\quad}}
                 &   & \go^0 &  \go^{-1} &  \go^{-2} & \go^{-3}\\
\hline
\mbox{(a)~data acquisition}  & 2^0    & 1.110            & 1.010               & 1.110        \\
 				             & 2^{-3} & 1.111            & 1.100               & 1.011       \\
\hline
\mbox{(b)~alignment} & 2^0  & 1.110            & 1.010               & 1.110        \\
 				     & 2^0  & 0.001            & 1.111               & 1.001    &  0.110    \\
\hline
\mbox{(c)~sum}       & 2^0  & 1.111            &11.001               &10.111    &  0.110\\
\hline
\mbox{(d)~redistribution} & 2^0  & 1.111             &                    &      &     \\
\hline
\mbox{}					  & 2^0  & 0.001  &1.001               &      & \\
\hline
\mbox{}					  & 2^0  &                   &0.001    &0.111 &  0.110\\
\hline
\mbox{(e)~sum}            & 2^0  & 10.000            &1.010               &0.111    &0.110\\
\hline
\mbox{(f)~normalization}  & 2^1 & 1.000            & 0.101               & 0.011&1.011     \\
\hline
\mbox{(g)~rounding}       & 2^1 & 1.000            & 0.101               & 0.100&     \\
\end{array}
\label{tab1}
$$
\end{table*}
First of all, the two numbers are stored in memory by distributing their digits along the powers $\go^0$, $\go^{-1}$ and $\go^{-2}$ (step (a)). Before summing the two numbers, an alignment is performed to make the two exponents equal (step (b)). Notice that shifting to the right the digits of the second number causes a redistribution of the digits along the three mantissas. Step (c) performs a polynomial-like sum of the two numbers. The contribution of each term has to be consistently redistributed (step (d)), in order to take into account possible carry bits,  and the three mantissas accordingly updated (step (e)). Steps (f) and (g) conclude the computation by normalizing and rounding the result.           

\paragraph{Subtraction.} As usual, floating-point subtraction bet\-ween two numbers sharing the same sign is performed  by inverting the sign bit of the second number, converting to 2's complement its mantissa, and then performing the addition as outlined above. It is well-known that subtracting two close numbers may lead to cancellation issues. We consider an example where the accuracy may be dynamically changed in order to overcome ill-conditioning issues. We assume to work with the arithmetic resulting  by setting  $t=7$  and $T=3$ (four grossdigits each consisting of one byte) with truncation. It turns out that, for a floating-point number $X$ representing an input real number $x$, its section  $X^{(0)}$ may be interpreted as the single precision representation of $x$,  while $X^{(1)}$, $X^{(2)}$ and $X^{(3)}\equiv X$ are its double, triple and quadruple precision approximations respectively. 
Loss of accuracy, resulting from a subtraction between two numbers having the same sign, will be detected during the normalization phase, when it requires shifting the mantissa by a large number of bits.

Consider the simple problem of evaluating the function $f(x,y,z)=x+y+z$ that computes the sum of three  real numbers, and assume that the user requires a simple precision accuracy in the result. In the examples below, we discuss three different situations.
\begin{example}
\label{subtraction_ex1}
The three real numbers 
$$
\begin{array}{rclcl}
x &=&  2^{-1} &\cdot& 1.0001100000010111111001001110110\cdots, \\
y &=&  2^0 &\cdot&    1.0010101010110010110101001101011\cdots,\\
z &=&  2^0 &\cdot&    1.1011011010111011011011010111001\cdots,
\end{array}
$$
are represented on the Infinity Computer as
$$
\begin{array}{rcll}
X &=& 2^{-1} \cdot (&\go^0 1.0001100 \,+\,\go^{-1} 0.0001011\\  
  &&& +\,\go^{-2} 0.1101010 \,+\, \go^{-3} 0.1101011),  \\
\end{array}
$$
$$
\begin{array}{rcll}
Y &=& 2^0    \cdot (&\go^0 1.0010101 \,+\, \go^{-1} 0.1011001 \\ 
&&& +\, \go^{-2} 1.0110110 \,+\, \go^{-3} 0.1101011),  \\
\end{array}
$$
$$
\begin{array}{rcll}
Z &=& 2^{0}  \cdot (&\go^0 1.1011011 +\, \go^{-1} 0.1011101\\
&&&+\, \go^{-2} 1.0110110 \,+\, \go^{-3} 1.0111001).
\end{array}
$$
Since we are adding positive numbers, no control on the accuracy is needed here, and the result is yielded as 
$$
f(X,Y,Z)\approx X^{(0)}+Y^{(0)}+Z^{(0)}= 2^1 \cdot 1.1011011,
$$ 
with a relative error $E^{(0)}\approx 1.1 \cdot 2^{-10}$, as is expected in simple precision.
\end{example}
\begin{example}
\label{subtraction_ex2}
Given the three real numbers defined in the previous example, we want now to evaluate $f(x,y,-z)$ again requiring  an eight-bit accuracy in the result. Table \ref{tab2} shows the sequence of steps performed to achieve the desired result.
\begin{table*}
\caption{Avoiding cancellation issues when evaluating the function $f(x,y,z)=x+y+z$ for the input data in Example \ref{subtraction_ex2}.}
$$
\begin{array}{ll@{\quad}|c@{\quad}|l@{\quad}}
& \multicolumn{1}{c|}{\mbox{steps}}    &  \multicolumn{1}{c|}{\mbox{error}}  &  \multicolumn{1}{c}{\mbox{action}}      \\          
\hline
\mbox{(a) } &S^{(0)} :=X^{(0)}+Y^{(0)}=2^0 \cdot \go^0 1.1011011 & <2^{-8} &   \mbox{ accept the result} \\
\mbox{    } &S^{(0)}-Z^{(0)}= 0 & 1 &   \mbox{ improve the accuracy} \\
\hline
\mbox{(b) } &S^{(1)} :=X^{(1)}+Y^{(1)}= 2^0 \cdot (\go^0 1.1011011 + \go^{-1}0.1011110)  & <2^{-16} &   \mbox{ accept the result} \\
\mbox{    } &S^{(1)}-Z^{(1)}= 2^{0} \cdot  \go^{-1}0.0000001 = 2^{-15} \cdot  \go^{0}1.0000000 & >2^{-2} &   \mbox{ improve the accuracy} \\
\hline
\mbox{(c) } &S^{(2)} :=X^{(2)}+Y^{(2)}= 2^0 \cdot (\go^0 1.1011011 + \go^{-1}0.1011111 + \go^{-2} 01100011)  & <2^{-16} &   \mbox{ accept the result} \\ 
\mbox{    } &
\begin{array}{rcl}
S^{(2)}-Z^{(2)}&=& 2^{0} \cdot  (\go^{-1}0.0000001  + \go^{-2} 1.0101101) \\
&=& 2^{-15} \cdot  (\go^{0}1.1010110 + \go^{-1}1.000000) 
\end{array}
& <2^{-12} &   \mbox{ final result} \\
\end{array}
$$
\label{tab2}
\end{table*}
The computation in simple precision, as in the previous example, is described in step (a): it leads to a clear cancellation phenomenon and, once detected, the accuracy is improved by letting the $\go^{-1}$ terms enter into play (step (b)). However,  the relative error remains higher than the prescribed tolerance, and accuracy needs to be improved by also considering the  $\go^{-2}$ terms. The computation is then repeated at step (c) and the  correct result is finally achieved. Notice that, in performing steps (b) and (c), one can evidently exploit the work already carried out in the previous step. The overall procedure thus requires 6 additions/sutractions of grossdigits, the same that would be needed by directly working with a 24-bit register which, for this case, is the minimum accuracy requirement to obtain eight correct bit in the result. This means that no extra effort is introduced during the steps. As a further remark, we stress again that a parallelization through the steps is also possible, even though we will not  discuss this issue here.     
\end{example}
\begin{example}
\label{subtraction_ex3}
We want to evaluate $f(x,-y,-z)$ requiring an eight-bit accurate result, now choosing 
$$
\begin{array}{rclcl}
x &=&  2^0    &\cdot&    1.0010101101010111111001001110110\cdots, \\
y &=&  2^0    &\cdot&    1.0010100010110010110101001101011\cdots,\\
z &=&  2^{-7} &\cdot&    1.0101000011011110010010110001010\cdots.
\end{array}
$$
Table \ref{tab3} shows the sequence of steps performed to achieve the desired result for this case.
\begin{table*}
\caption{Avoiding cancellation issues when evaluating the function $f(x,y,z)=x+y+z$ for the input data in Example \ref{subtraction_ex3}.}
$$
\begin{array}{ll@{\quad}|c@{\quad}|l@{\quad}}
& \multicolumn{1}{c|}{\mbox{steps}}    &  \multicolumn{1}{c|}{\mbox{error}}  &  \multicolumn{1}{c}{\mbox{action}}      \\          
\hline
\mbox{(a) } &S^{(0)} :=X^{(0)}-Y^{(0)}=2^{-7} \cdot \go^0 1.0000000 & \approx 2^{-2} &   \mbox{ improve the accuracy} \\
\hline
\mbox{(b) } &S^{(1)} :=X^{(1)}-Y^{(1)}= 2^{-7}  \cdot \go^0 1.0101001   & <2^{-9} &   \mbox{ accept the result} \\
\mbox{    } &S^{(1)}-Z^{(0)}= 2^{-14} \cdot  \go^{0}1.0000000  & >2^{-1} &   \mbox{ improve the accuracy} \\
\hline
\mbox{(c) } &S^{(2)} :=X^{(2)}-Y^{(2)}= 2^{-7} \cdot (\go^0 1.0101001 + \go^{-1}0.1000100)  & <2^{-20} &   \mbox{ accept the result} \\ 
\mbox{    } &
S^{(2)}-Z^{(1)}= 2^{-7} \cdot  \go^{-1}1.1010101  = 2^{-15} \cdot \go^{0}1.1010101  
& <2^{-10} &   \mbox{ final result} \\
\end{array}
$$
\label{tab3}
\end{table*}
When working in simple precision, an accuracy improvement is already needed when subtracting the first two terms  $X^{(0)}$ and $Y^{(0)}$ and, consequently, step (a) is stopped. At step (b), the difference $x-y$ is evaluated in double precision which, on balance, assures an eight-bit accuracy in the result. However, a new cancellation issue emerges when  subtracting  $Z^{(0)}$ from $X^{(1)}-Y^{(1)}$, suggesting that the two terms need to be represented more accurately. This is done in step (c), evaluating $x-y$ in triple precision and representing $z$ in double precision. The overall procedure requires 5 additions/sutractions of grossdigits. This example, compared with the previous one, reveals the coexistence of variables combined with different precisions.  
\end{example}

Summarizing the three examples above, we observe how the accuracy of representation and combination of variables may be dynamically changed, in order to overcome possible loss of significant figures in the result when evaluating a function. Of course, for this strategy to work, it is necessary that the input data are stored with high precision and a technique to detect the loss of accuracy be available. In Section \ref{sec:5} we will illustrate this procedure applied to the accurate determination of zeros of functions (a further example may be found in  \cite{AmBrIaMa20}).   

Concerning the computational complexity, it should be noticed that Example \ref{subtraction_ex1} reflects the normal situation where the use of the standard precision is enough to produce a correct result, while Examples \ref{subtraction_ex2} and \ref{subtraction_ex3} highlight less frequent events. 

\paragraph{Multiplication.} 
Set $t=3$  and $T=2$ (three grossdigits each with four significant digits).   
Consider the product of the two floating-point normalized numbers
$$
\begin{array}{rcl}
X &=&  2^0    \cdot 1.01101111100,  \\
Y &=&  2^{0}  \cdot 1.10111111101. 
\end{array}
$$
Table \ref{tab4} summarizes the procedure by a sequence of commented steps.
\begin{table*}
\caption{Scheme of the multiplication of two  floating-point numbers.}
$$
\begin{array}{l@{\quad}|l@{\quad}|l@{\quad}|l@{\quad}|l@{\quad}|l@{\quad}|l@{\quad}|l@{\quad}|l@{\quad}}
                 &   & \go^1   & \go^0 &  \go^{-1} &  \go^{-2} &  \go^{-3} &  \go^{-4} & \go^{-5}\\
\hline
\mbox{(a)~data acquisition} & 2^0 &  & 1.011  & 0.111  & 1.100  &&&  \\
 				            & 2^0 &  & 1.101  & 1.111  & 1.101  &&&  \\
\hline
\mbox{(b)~convolution product} & 2^0 & & 10.001111 & 100.000000 & 110.010100 & 100.001111 & 10.011100 & \\
\hline
\mbox{(c)~redistribution} & 2^0 & 0.001 & 0.001 & 1.110 &       &        &       &     \\
\hline
\mbox{}					  & 2^0 &       & 0.010 & 0.000 & 0.000 &        &       &     \\
\hline
\mbox{}					  & 2^0 &       &       & 0.011 & 0.010 & 1.000  &       &     \\
\hline
\mbox{}					  & 2^0 &       &       &       & 0.010 & 0.001  & 1.110 &     \\
\hline
\mbox{}					  & 2^0 &       &       &       &       & 0.001  & 0.011 & 1.000\\
\hline
\mbox{(d)~sum with redistribution} & 2^0 & 0.001 & 0.100 & 0.001 & 0.100 & 1.011 & 0.001 &  1.000  \\
\hline
\mbox{(e)~normalization}  & 2^1 &        & 1.010 & 0.000 & 1.010 & 0.101 & 1.000 & 1.100   \\
\hline
\mbox{(f)~rounding}  & 2^1 &        & 1.010 & 0.000 & 1.010 &       &       &         \\
\end{array}
\label{tab4}
$$
\end{table*}
After expanding the input data along the negative powers of $\go$ for data storage (step (a)), the convolution product described in (\ref{Xq*Yq}) is performed (step (b)). At step (c), the contribution of each term is redistributed, and a sum is then needed to update the mantissas (step (d)). Steps (e) and (f) conclude the computation by normalizing and rounding the result. Notice that step (e) may be carried out by applying the rules for the addition described in Table \ref{tab1}. Again, we stress that the terms in the convolution product, as well as in the subsequent sum, may be computed in parallel.
\paragraph{Division.} 
The division of two floating-point numbers $X$ and $Y$ has been switched to the multiplication of $X$ by the reciprocal of $Y$. This latter, in turn, is obtained with the aid of the Newton-Raphson method applied to find the zero of the function $f(Z) = 1/Z - Y$. Hereafter, without loss of generality, we assume $Y>0$. Starting from a suitable initial guess $Z_0$, the Newton iteration then reads
\begin{equation}
\label{newton}
Z_{k+1}=Z_k+Z_k(1-Y Z_k).
\end{equation}
The relative error 
$$
E_k:=\frac{1/Y-Z_k}{1/Y}=1-YZ_k
$$ 
satisfies
\begin{equation}
\label{err_newt}
E_{k+1}= 1-YZ_{k+1}=1-2YZ_k+(YZ_k)^2=E_k^2,
\end{equation}
which means that, as is expected in presence of simple zeros, the sequence $Z_k$ eventually converges quadratically to $1/Y$, and the number of correct figures doubles at each iteration. This feature makes the division procedure extremely efficient in our context, since the required accuracy may be easily increased to an arbitrary level. In order to obtain such a good convergence rate starting from the very beginning of the sequence, the numerator $X$ and denominator $Y$ are scaled by a suitable factor $\beta^s$ so that $\widehat Y:=\beta^s Y$ lies in the interval $[0.5,\,1]$. In the literature, the minmax linear polynomial approximation is often used to estimate the reciprocal of $\widehat Y$. The resulting initial guess is
$$
Z_0=\frac{48}{17}-\frac{32}{17}\widehat Y,
$$
which assures an initial error $E_0\le 1/17$. Taking into account the equality (\ref{err_newt}), the relative error at step $k$ decreases as 
$$
E_k = E_0^{2^k} \le \left( \frac{1}{17} \right)^{2^k},
$$
and consequently, assuming $\beta=2$, a $q$-bit accurate approximation is obtained by setting
$$
k = \left\lceil \log_2 \frac{q+1}{\log_2 17}\right\rceil,
$$
where $\lceil \cdot \rceil$ denotes the ceiling function. As an example, four iterations suffice to get an approximation with at least $32$ correct digits. 
Table \ref{tab5} shows the sequence generated from the scheme above applied to find, on the Infinity Computer, the reciprocal of the binary number $Y=(1010)_2$ ($1/10$ in decimal base), under the choice $t=3$ and $T=7$ (eight grossdigits each with four significant figures).    
\begin{table*}
\caption{Newton iteration to compute the reciprocal of $Y=2^0\cdot 1010$ on the Infinity Computer.}
$$
\begin{array}{c@{\quad}|@{~}c@{~}|l@{\quad}|l@{\quad}|l@{\quad}|l@{\quad}|l@{\quad}|l@{\quad}|l@{\quad}|l@{\quad}}
\mbox{sequence}      &  & \go^0   & \go^{-1}   & \go^{-2} &  \go^{-3} &  \go^{-4} &  \go^{-5} &  \go^{-6} & \go^{-7}\\[.02cm]
\hline
Z_0 & 2^0  & 1.010 & 0.000 & 0.000 & 0.000 & 0.000 & 0.000 & 0.000 & 0.000 \\[.02cm]
Z_1 & 2^0  & 1.101 & 0.100 & 1.101 & 0.100 & 1.101 & 0.100 & 1.101 & 0.100 \\[.02cm]
Z_2 & 2^0  & 1.100 & 1.100 & 0.111 & 1.100 & 0.010 & 1.011 & 1.101 & 1.100 \\[.02cm]
Z_3 & 2^0  & 1.100 & 1.100 & 1.100 & 1.100 & 1.010 & 1.101 & 0.000 & 1.110\\[.02cm]
Z_4 & 2^0  & 1.100 & 1.100 & 1.100 & 1.100 & 1.100 & 1.100 & 1.100 & 1.101  
\end{array}
\label{tab5}
$$
\end{table*}
\subsection{Implementation details}
\label{sec:4}
We have developed a Matlab prototype emulating the Infinity Computer environment interfaced with a module that performs the suitable carrying, normalization and rounding processes, needed by the identification of $\go$ and $\gob$ to ensure proper functioning of the resulting dynamic floating-point arithmetic. 

The emulator  represents input real numbers using a set of binary grossdigits, whose length and number are defined by the two input parameters $t$ and $T$.  This latter parameter is  used to define the maximum available accuracy for storing variables. In accord with formulae such as (\ref{Xq+Yq}) and (\ref{Xq*Yq}), the actual accuracy used to execute a single operation will depend on the accuracy of the two operands but cannot exceed $T$.

At the moment, the emulator implements the four basic operations following the strategies described above, plus some simple functions. The vectorization issue, to speed-up the execution time associated with each floating-point operation, has not yet been addressed, so that  all operations between grossdigits are executed sequentially. 

All computations reported in the present paper, including the results presented in the next section,  have been carried out  on an Intel i5 quad-core computer with 16GB of memory, running Matlab R2019b.

\section{A numerical illustration}
\label{sec:5}
As an application highlighting the potentialities of the dynamic precision arithmetic introduced above, we consider the problem of determining accurate approximations of the zeros of a function $f:[a,b]\rightarrow \RR$,  in the case where this problem suffers from ill-conditioning issues. 

The finite arithmetic representation of the function $f$ introduces perturbation terms of different nature: analytical errors, errors in the coefficients or parameters involved in the definition of the function, or roundoff errors introduced during its evaluation. 

From a theoretical point of view, these sources of errors  may be accounted for by introducing a perturbation function $g(x)$ and analyzing its effects on the zeros of the perturbed  function  $\tilde f(x):= f(x)+\eps g(x)$ where the factor $\eps$ has the size of the unit roundoff. Under regularity assumptions on $f$, if $\alpha\in (a,b)$ is a zero with multiplicity $d>0$, it turns out that $\tilde f(x)$ admits a perturbed zero $\alpha+\delta \alpha $, with the perturbing term $\delta \alpha$ satisfying, in first approximation,
\begin{equation}
\label{delta_alpha}
|\delta \alpha| \approx \eps^{1/d} d! \left| \frac{g(\alpha)}{f^{(d)}(\alpha)} \right|^{1/d}.
\end{equation}
As an example, consider the polynomial  
\begin{equation}
\label{pol1}
p(x)=x^5 -5x^4+10x^3 -10x^2 +5x -1
\end{equation}
that admits $\alpha =1$ as unique root with multiplicity $d=5$ (indeed $p(x)=(x-1)^5$). 
For this problem, from formula (\ref{delta_alpha}) we get   
\begin{equation}
\label{delta_alpha_pol1}
\left| \frac{\delta \alpha}{\alpha}\right| = |\delta \alpha| \approx \eps^{1/d} |g(1)|^{1/d}.
\end{equation}
Working with 64-bit IEEE arithmetic, i.e. with a roundoff unit $u=2^{-53}$, we expect a breakdown of the relative error proportional to $u^{1/5} \approx 6.4 \cdot 10^{-4}>0.5 \cdot 10^{-3}$, so that, assuming $|g(1)|^{1/d}\approx 1$, the approximation of the zero $\alpha$ only contains $3\div4$ correct figures. 

\begin{figure}
\hspace*{-.3cm}  \includegraphics[width=0.52\textwidth]{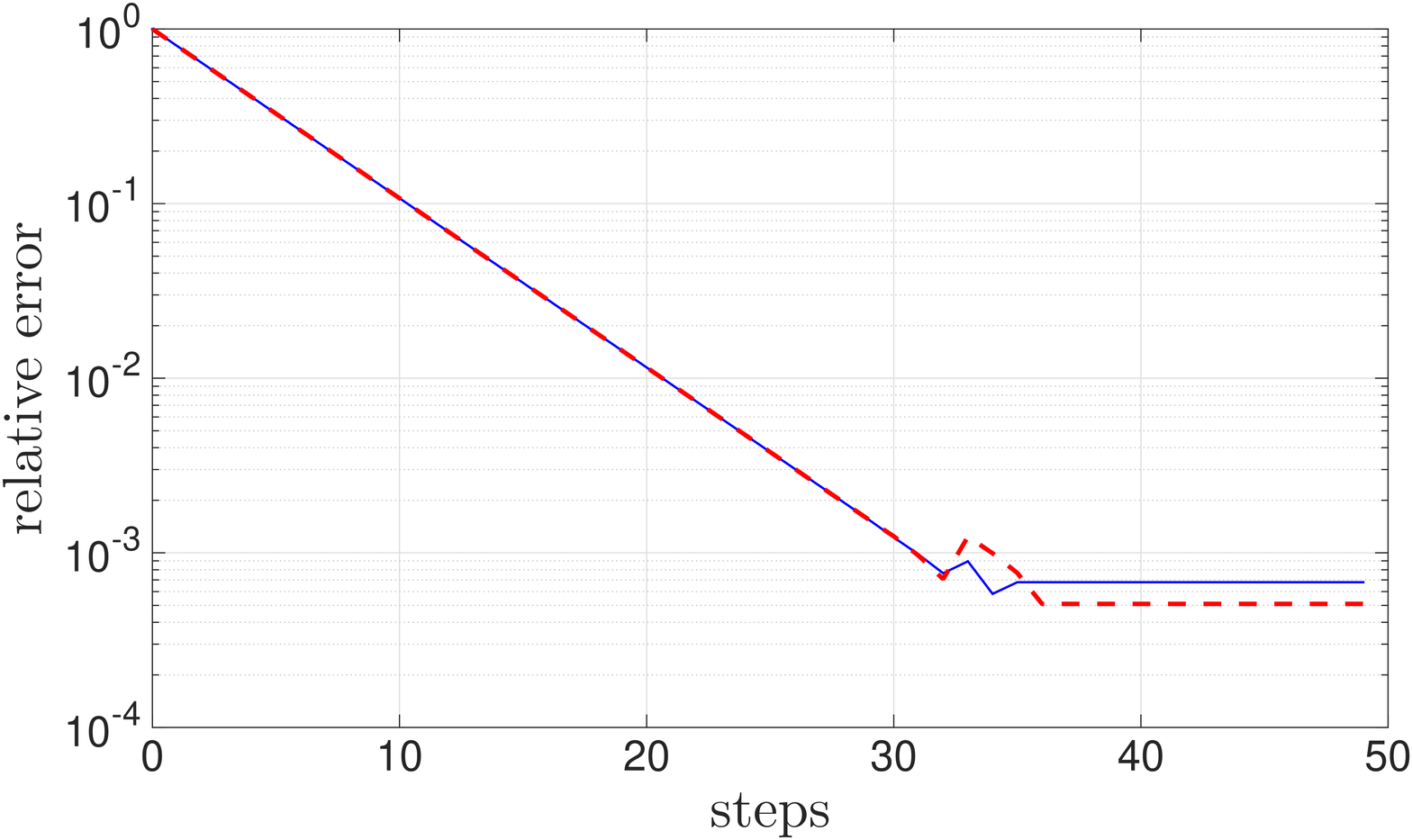}
\caption{Relative error related to the sequence of approximations generated by the Newton method applied to the polynomial (\ref{pol1}). Solid line: implementation on the Infinity Computer with $t=52$ and $T=0$. Dashed line: implementation in Matlab double precision arithmetic.}
\label{fig1}       
\end{figure}

This is confirmed by the two plots in Figure  \ref{fig1}. They display the relative error of the approximations to $\alpha$ generated by applying the Newton method to the problem $p(x)=0$, choosing $x_0=2$ as initial guess:
\begin{equation}
\label{newton_pol1}
x_{k+1}=x_k-\frac{p(x_k)}{p'(x_k)}.
\end{equation}
The solid line refers to the implementation of the iteration on the Infinity Computer using $t=52$ and $T=0$. This choice mimics the default double precision arithmetic in Matlab, which uses a register of 64 bit to store a normalized binary number, 52 bit being dedicated to the (fractional part of the) mantissa. As a matter of fact, the dashed line, coming out from the implementation of the scheme using the standard Matlab arithmetic, precisely overlap with the solid line as long as the error decreases, while the two lines slightly depart from each other when they reach the saturation level right below $10^{-3}$, namely starting from step $32$.   

\begin{figure}
\hspace*{-.3cm}  \includegraphics[width=0.52\textwidth]{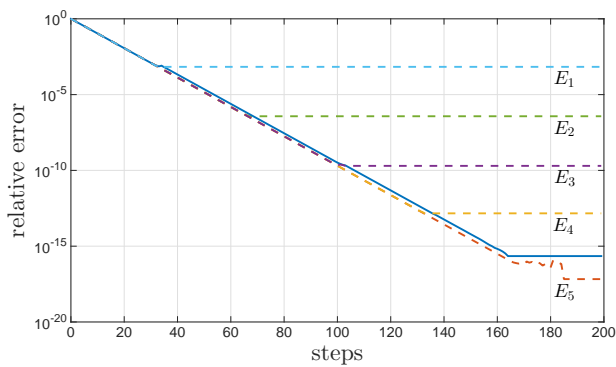}
\caption{Relative error corresponding to the sequence of approximations generated by the Newton method applied to the polynomial (\ref{pol1}) on the Infinity Computer. Solid line: dynamic precision implementation. Dashed line: fixed precision implementation, for different accuracies.}
\label{fig2}       
\end{figure}

We want now to improve the accuracy of the approximation to the zero $\alpha=1$ of (\ref{pol1}) by exploiting the new computational platform. 
Hereafter, the $53$-bit  precision used above will be referred to as \textit{single precision}.  The dashed lines in Figure \ref{fig2} show the relative error reduction when the Newton method is implemented on the Infinity Computer by working with multiple fixed precision. From top to bottom, we can see the five saturation levels corresponding to the stagnation of the error at  
$E_1\approx 6.8\cdot 10^{-4}$ in single precision,  $E_2\approx 3.7\cdot 10^{-7}$ in double precision, $E_3\approx 2.0 \cdot 10^{-10}$ in triple precision, $E_4\approx 1.5 \cdot 10^{-13}$ in quadruple precision, and $E_5\approx 6.7 \cdot 10^{-18}$ in quintuple precision. These saturation values are consistently  predicted by formula (\ref{delta_alpha_pol1}), after replacing $\eps$ with $2^{-53k}$, for $k=1,\dots,5$. 

Now suppose we want $53$ correct binary digits in the approximation (i.e., about $15\div 16$ correct decimal digits). From the discussion above, it turns out that we have to activate the quintuple precision, thus setting  $t=52$ and $T=4$ (five grossdigits, each consisting of a 53-bit register). However, the computational effort may be significantly reduced if we increase the accuracy by involving new negative grosspowers only when they are really needed. In a dynamic usage of the accuracy, starting from $x_0$, we can initially activate the single precision mode until we reach the first saturation level and, thereafter, switch to double precision until the second saturation level is reached, and so forth until we get the desired accuracy in the approximation. Denoting by 
$$
\mbox{err}(k) =\left| \frac{x_{k}-x_{k-1}}{x_k} \right|
$$
the estimated error at step $k$, and by \texttt{prec} the current precision, initially set equal to $1$,  the points where an increase of the accuracy is needed may be  automatically detected by employing a simple control scheme such as
\begin{verbatim}
if err(k)>=s*err(k-1) and prec <=T
   prec=prec+1
end   
\end{verbatim}   
where $s\le 1$ is a positive safety factor that we have set equal to $1$. The solid line in  Figure \ref{fig2} shows the corresponding reduction of the error and we can see that the change of precision scheme described above works quite well for this example, since all saturation levels are correctly detected and overcome. At step $162$ the error reaches its minimum value of $2.2 \cdot 10^{-16}$ and the iteration could be stopped by the standard criterion $\mbox{err}(k)<10^{-15}$ even though, for clarity, we have generated additional points to reveal the last saturation level corresponding to \texttt{prec}$=T+1=5$.     

Now, let us compare the computational cost of the dynamic implementation versus the fixed quintuple precision one, considering that to reach the highest precision each mode requires $162$ Newton iterations (see Figure \ref{fig2}). On the basis of the formula reported right below (\ref{Xq*Yq}),  the dynamic implementation would take about $2.4 \cdot 10^3$ grossdigits multiplications while the fixed quintuple precision implementation requires $2.0\cdot 10^4$ grossdigits multiplications.\footnote{For simplicity, we do not consider additions/subtractions in the computation, since their contribution would not alter the final result.} It follows that the former mode would reduce the execution times of a factor at least eight with respect to the latter. Actually, it does much better: the dynamic usage of variables and operations, understood as the ability of handling variables with different accuracy and executing operations on them, makes the resulting arithmetic definitely much more efficient than what emerged from the comparison above. 
\begin{table*}
\caption{The Horner method for evaluating $p(x)$ in (\ref{pol1}) at $x=2^0\cdot 1.0000000000000000000000000000000000000000000001000010$.}
$$
\begin{array}{l@{\quad}|c|l|l@{\quad}}
{p=1}           & 2^0  & \go^0  & 1.0000000000000000000000000000000000000000000000000000 \\[.05cm]
\hline 
{p=p\cdot x-5}  & -2^1 & \go^0  & 1.1111111111111111111111111111111111111111111111011111 \\[.05cm]
\hline
{p=p\cdot x+10} & 2^2 & \go^0  & 1.0111111111111111111111111111111111111111111111001110 \\[.05cm]
                &     &\go^{-1}& 1.0000000000000000000000000000000000000000100010000010 \\[.05cm]
\hline
{p=p\cdot x-10} &-2^1 & \go^0  & 1.1111111111111111111111111111111111111111111110011101 \\[.05cm]
                &     &\go^{-1}& 0.0000000000000000000000000000000000000010001000000111 \\[.05cm]
                &     &\go^{-2}& 1.1111111111111111111111111111111101110011100111110000 \\[.05cm]
\hline
{p=p\cdot x+5} &2^{-1}& \go^0  & 1.1111111111111111111111111111111111111111111101111100 \\[.05cm]
                &     &\go^{-1}& 0.0000000000000000000000000000000000000100010000001111 \\[.05cm]
                &     &\go^{-2}& 1.1111111111111111111111111111110111001110011111000000 \\[.05cm]
                &     &\go^{-3}& 0.0000000000000000000000010010000110001000000100000000 \\[.05cm]
\hline
{q=p\cdot x}	&2^{0}& \go^0  & 1.0000000000000000000000000000000000000000000000000000 \\[.05cm]
                &     &\go^{-1}& 0.0000000000000000000000000000000000000000000000000000 \\[.05cm]
                &     &\go^{-2}& 0.0000000000000000000000000000000000000000000000000000 \\[.05cm]
                &     &\go^{-3}& 0.0000000000000000000000000000000000000000000000000000 \\[.05cm]
                &     &\go^{-4}& 0.0000000000000000010010101010010100010100001000000000 \\[.25cm]
p=p-1       & 2^{-230}&\go^0   & 1.0010101010010100010100001000000000000000000000000000
\end{array}
\label{tab6}
$$
\end{table*}


In carrying out the computation above, for the dynamic precision mode we have assumed that all floating-point operations were executed with the current selected precision. For example, under this assumption, the computational effort per step of the two modes would become equivalent starting from step $139$ onwards since, at that step, the dynamic mode activates the quintuple precision to overcome the threshold level $E_4$ in Figure \ref{fig2}. 

There is, however, one fundamental aspect that we have not yet considered. In fact, to overcome the ill-conditioning of the problem, the higher precision is only needed during the evaluation of $p(x_k)$ and $p'(x_k)$ in (\ref{newton_pol1}), while the single $53$-bit precision is enough to handle the sequence $x_k$. In other words, to minimize the overall computational effort, we may improve the accuracy only in the part of the code that implements the Horner rule to evaluate the polynomial $p(x)$ and its derivative. 

Interestingly, we have not to instruct the Infinity Computer to switch between single and quintuple precision: all is done automatically and naturally and, more importantly, even during the evaluation of $p(x_k)$ and $p'(x_k)$, the transition from single to quintuple precision is gradual, in that all the intermediate precisions are actually involved only when really needed, which makes the whole machinery much more efficient. 

To better elucidate this aspect, we illustrate the sequence produced by the Horner rule to evaluate $p(x_k)$ at step $k=145$, where the quintuple precision is activated. The first column in Table \ref{tab6} reports the five steps of the Horner method applied to evaluate the polynomial $p(x)$ in (\ref{pol1}) at the floating-point single precision number $x=x_{145}$ (its value is in the caption of the table). The variable $p$ is initializated with the leading coefficient of the polynomial, but is allowed to store five grossdigits, each $53$-bit long, to host  floating-point numbers up to quintuple precision. From the table we see that, as the iteration scheme proceeds, new negative grosspowers appear in the values taken by  the variable $p$. More precisely, at step $k$ the variable $p$ stores a $k$-fold precision floating-point number, for $k=1,\dots,5$. 

The increase in the precision of one unit at each step evidently arises from the product $p\cdot x$, since $x$ remains a single-precision variable and no rounding occurs.  Let us better examine what happens at the last step. The product $p\cdot x$ generates a quintuple-precision number whose expansion along negative grosspowers  matches the number $1$ up to $\go^{-3}$. Consequently, the last operation $p-1$ only contains significant digits in the coefficient of $\go^{-4}$ so that, after normaliziation, $p$ will store again a single-precision number that can be consistently combined inside formula (\ref{newton_pol1}). 

In conclusion, the Horner procedure, though being enabled to operate in quintuple precision, actually involves lower precision numbers, except at the very last step. The five steps reported in Table \ref{tab6} require $15$ multiplications of grossdigits, with a clear saving of time, if we consider that the fixed quintuple-precision mode would require $125$ multiplications of grossdigits. Comparing the execution times in Matlab over $162$ steps, we found out that the dynamic-precision implementation is about $1.75$ times slower than the single-precision implementation (which however stagnates at level $E_1$) and about $19$ times faster than the quintuple precision mode, thus confirming the expected efficiency.

\section{Conclusions}
\label{sec:6}
We have proposed a variable precision floating-point arithmetic able to simultaneously storing numbers and execute operations with different accuracies. This feature allows one to dynamically change the accuracy during the execution of a code, in order to prevent inherent ill-conditioning issues associated with a given problem. In this context, the Infinity Computer has been recognized as a natural computational environment that can easily host such an arithmetic. The assumption that makes this paradigm work is the identification of the two symbols $\go$ and $\gob$. The latter, defined as $\gob=\beta^{t+1}$, is evidently a finite quantity for our numeral system but, in many respects, its reciprocal behaves as an infinitesimal-like entity in the numeral system induced by a floating-point arithmetic operating with $t+1$ significant figures. In the same spirit of the Infinity Computer, it turns out that negative powers of $\gob$ may be used as ``lenses'' to increase and decrease the accuracy when needed. An emulator of this dynamic precision floating-point arithmetic has been developed in Matlab, and an application to the accurate solution of (possibly ill-conditioned) scalar nonlinear equations has been discussed.

\begin{acknowledgements}
This work was funded by the INdAM-GNCS 2018 Research Project ``Numerical methods in optimization and ODEs'' (the authors are members of the INdAM Research group GNCS).
\end{acknowledgements}

%
%

\bibliographystyle{spmpsci}      
\bibliography{dynprec}   

%
%


\end{document}